\documentclass[12pt]{amsart}
\usepackage{amssymb}
\usepackage{amsfonts}

\newtheorem{theorem}{Theorem}

\theoremstyle{definition}

\theoremstyle{remark}

\begin{document}

\author{E. Liflyand}

\title {On Fourier re-expansions}



\address{Department of Mathematics, Bar-Ilan University, 52900 Ramat-Gan, Israel}
\email{liflyand@math.biu.ac.il}

\begin{abstract} We study an extension to Fourier transforms of
the old problem on absolute convergence of the re-expansion in the
sine (cosine) Fourier series of an absolutely convergent cosine
(sine) Fourier series. The results are obtained by revealing certain
relations between the Fourier transforms and their Hilbert
transforms.
\end{abstract}

\maketitle

\section{Introduction}

In 50-s (see, e.g., \cite{IT} or in more detail \cite[Chapters II
and VI]{Kah}), the following problem in Fourier Analysis attracted
much attention:

{\it Let $\{a_k\}_{k=0}^\infty$ be the sequence of the Fourier
coefficients of the absolutely convergent sine (cosine) Fourier
series of a function $f:\mathbb T=[-\pi,\pi)\to \mathbb C,$ that is
$\sum |a_k|<\infty.$ Under which conditions on $\{a_k\}$ the
re-expansion of $f(t)$ ($f(t)-f(0)$, respectively) in the cosine
(sine) Fourier series will also be absolutely convergent?}

The obtained condition is quite simple and is the same in both
cases:

\begin{eqnarray}\label{condser}\sum\limits_{k=1}^\infty
|a_k|\ln(k+1)<\infty. \end{eqnarray}

In this paper we study a similar problem for Fourier transforms
defined on $\mathbb R_+=[0,\infty).$ Let

\begin{eqnarray*} F_c(x)=\int_0^\infty f(t)\cos xt\,dt\end{eqnarray*}
be the cosine Fourier transform of $f$ and

\begin{eqnarray*}F_s(x)=\int_0^\infty f(x)\sin xt\,dt\end{eqnarray*}
be the sine Fourier transform of $f$, each understood in certain
sense.

Let

\begin{eqnarray}\label{ac}\int_0^\infty|F_c(x)|\,dx<\infty,\end{eqnarray}
and hence

\begin{eqnarray}\label{ac1}f(t)=\frac{1}{\pi}\int_0^\infty F_c(x)\cos tx\,dx,\end{eqnarray}
or, alternatively,

\begin{eqnarray}\label{as}\int_0^\infty|F_s(x)|\,dx<\infty\end{eqnarray}
and hence

\begin{eqnarray}\label{as1}f(t)=\frac{1}{\pi}\int_0^\infty F_s(x)\sin tx\,dx.\end{eqnarray}

1) {\it Under which (additional) conditions on $F_c$ we get}
(\ref{as}),

\noindent or, in the alternative case,

2) {\it under which (additional) conditions on $F_s$ we get}
(\ref{ac})?

\noindent The answer is given by the following theorem. To formulate
it we turn to the Hilbert transform of an integrable function $g$

\begin{eqnarray}\label{dht}\mathcal{H}g(x)=\frac{1}{\pi}\int_\mathbb{R}\frac{g(t)}{x-t}\,dt,\end{eqnarray}
where the integral is understood in the improper (principal value)
sense, as $\lim\limits_{\delta\to0+}\int_{|t-x|>\delta}.$ It is not
necessarily integrable, and when it is, we say that $g$ is in the
(real) Hardy space $H^1(\mathbb R).$ If $g\in H^1(\mathbb R)$, then

\begin{eqnarray}\label{vm}\int_{\mathbb R}  g(t)\,dt=0.\end{eqnarray}
It was apparently first mentioned in \cite{kober}.

\begin{theorem}\label{ans} In order than the re-expansion $F_s$
of $f$ with the integrable cosine Fourier transform $F_c$ be
integrable, it is necessary and sufficient that its Hilbert
transform ${\mathcal H}F_c(x)$ be integrable.

Similarly, in order than the re-expansion $F_c$ of $f$ with the
integrable sine Fourier transform $F_s$ be integrable, it is
necessary and sufficient that its Hilbert transform ${\mathcal
H}F_s(x)$ be integrable.
\end{theorem}

\bigskip

\section{Proof of Theorem \ref{ans}}

Let (\ref{ac}) holds true. Then we can rewrite

\begin{eqnarray}\label{rw1} F_s(x)=\int_0^\infty\,\biggl[\frac{1}{\pi}
\int_0^\infty F_c(u)\cos tu\,du\biggr]\,\sin xt\,dt.
\end{eqnarray}
The right-hand side can be understood in the $(C,1)$ sense as

\begin{eqnarray*}\frac{1}{\pi}\lim\limits_{N\to\infty}\int_0^N\,
(1-\frac{t}{N})\int_0^\infty F_c(u)\cos tu\,du\,\sin xt\,dt.
\end{eqnarray*}
In virtue of (\ref{ac}) we can change the order of integration:

\begin{eqnarray*}&\quad&\frac{1}{\pi}\lim\limits_{N\to\infty}\int_0^\infty F_c(u)
\int_0^N\, (1-\frac{t}{N})\cos tu\,\sin xt\,dt\,du\\
&=&\frac{1}{\pi}\lim\limits_{N\to\infty}\int_0^\infty F_c(u)\,
\frac12 \int_0^N\, (1-\frac{t}{N})\,[\sin(u+x)t-\sin(u-x)t]\,dt\,du.
\end{eqnarray*}
We now need the next simple formula

\begin{eqnarray}\label{ces}\int_0^N(1-\frac{t}{N})\sin At\,dt
=\frac{1}{A}-\frac{\sin NA}{NA^2}.             \end{eqnarray}
Applying it yields

\begin{eqnarray}\label{rw2} F_s(x)&=&\frac{1}{\pi}\lim\limits_{N\to\infty}\int_0^\infty F_c(u)\,
\biggl[\frac{1}{u+x}-\frac{\sin(u+x)N}{N(u+x)^2}\biggr]\,du\nonumber\\
&-&\frac{1}{\pi}\lim\limits_{N\to\infty}\int_0^\infty F_c(u)\,
\biggl[\frac{1}{u-x}-\frac{\sin(u-x)N}{N(u-x)^2}\biggr]\,du=I_1+I_2.
\end{eqnarray}

Let us begin with $I_2$. Substituting $u-x=t,$ we obtain

\begin{eqnarray*}I_2=-\frac{1}{\pi}\lim\limits_{N\to\infty}\int_{-x}^\infty F_c(x+t)\,
\biggl[\frac{1}{t}-\frac{\sin Nt}{Nt^2}\biggr]\,dt.\end{eqnarray*}
For $I_1,$ we first substitute $u=-v.$ Thus

\begin{eqnarray*}I_1&=&\frac{1}{\pi}\lim\limits_{N\to\infty}\int_{-\infty}^0 F_c(-v)\,
\biggl[\frac{1}{-v+x}-\frac{\sin(-v+x)N}{N(x-v)^2}\biggr]\,dv\\
&=&-\frac{1}{\pi}\lim\limits_{N\to\infty}\int_{-\infty}^0 F_c(-v)\,
\biggl[\frac{1}{v-x}-\frac{\sin(v-x)N}{N(v-x)^2}\biggr]\,dv\\
&=&-\frac{1}{\pi}\lim\limits_{N\to\infty}\int_{-\infty}^0 F_c(v)\,
\biggl[\frac{1}{v-x}-\frac{\sin(v-x)N}{N(v-x)^2}\biggr]\,dv.\end{eqnarray*}
The last equality follows from the evenness of $F_c.$ Substituting
$v-x=t,$ we obtain

\begin{eqnarray*}I_1=-\frac{1}{\pi}\lim\limits_{N\to\infty}\int_{-\infty}^{-x} F_c(x+t)\,
\biggl[\frac{1}{t}-\frac{\sin Nt}{Nt^2}\biggr]\,dt.\end{eqnarray*}
Therefore,

\begin{eqnarray*}F_s(x)=-\frac{1}{\pi}\lim\limits_{N\to\infty}\int_{-\infty}^\infty F_c(x+t)
\biggl[\frac{1}{t}-\frac{\sin Nt}{Nt^2}\biggr]\,dt.\end{eqnarray*}
We are now in a position to apply the following result (see
\cite[Vol.II, Ch.XVI, Th. 1.22]{Zg}; even more general result can be
found in \cite[Th.107]{Ti1}).

{\bf Theorem A.} {\it If $\displaystyle\frac{|f(t)|}{1+|t|}$ is
integrable on $\mathbb R$, then the $(C,1)$ means

\begin{eqnarray*}-\frac{1}{\pi}\int_{-\infty}^\infty f(x+t)
\biggl[\frac{1}{t}-\frac{\sin Nt}{Nt^2}\biggr]\,dt\end{eqnarray*}
converge to the Hilbert transform $\mathcal{H}f(x)$ almost
everywhere as} $N\to\infty$.

It follows from Theorem A that for almost all $x$)

\begin{eqnarray}\label{cosi}F_s(x)={\mathcal H}F_c(x).\end{eqnarray}
We remark that any integrable function satisfies the assumption of
Theorem A.

Now, let (\ref{as}) holds true. Then we can rewrite

\begin{eqnarray}\label{rw3} F_c(x)=\int_0^\infty\,\biggl[\frac{1}{\pi}
\int_0^\infty F_s(u)\sin tu\,du\biggr]\,\cos xt\,dt.
\end{eqnarray}
The right-hand side can be understood in the $(C,1)$ sense as

\begin{eqnarray*}\frac{1}{\pi}\lim\limits_{N\to\infty}\int_0^N\,
(1-\frac{t}{N})\int_0^\infty F_s(u)\sin tu\,du\,\cos xt\,dt.
\end{eqnarray*}
In virtue of (\ref{as}) we can change the order of integration:

\begin{eqnarray*}&\quad&\frac{1}{\pi}\lim\limits_{N\to\infty}\int_0^\infty F_s(u)
\int_0^N\, (1-\frac{t}{N})\sin tu\,\cos xt\,dt\,du\\
&=&\frac{1}{\pi}\lim\limits_{N\to\infty}\int_0^\infty F_s(u)\,
\frac12 \int_0^N\, (1-\frac{t}{N})\,[\sin(u+x)t+\sin(u-x)t]\,dt\,du.
\end{eqnarray*}
Applying (ref{ces}), we get

\begin{eqnarray}\label{rw4} F_s(x)&=&\frac{1}{\pi}\lim\limits_{N\to\infty}\int_0^\infty F_s(u)\,
\biggl[\frac{1}{u+x}-\frac{\sin(u+x)N}{N(u+x)^2}\biggr]\,du\nonumber\\
&+&\frac{1}{\pi}\lim\limits_{N\to\infty}\int_0^\infty F_s(u)\,
\biggl[\frac{1}{u-x}-\frac{\sin(u-x)N}{N(u-x)^2}\biggr]\,du=J_1+J_2.
\end{eqnarray}

Let us begin with $J_2$. Substituting $u-x=t,$ we obtain

\begin{eqnarray*}J_2=\frac{1}{\pi}\lim\limits_{N\to\infty}\int_{-x}^\infty F_s(x+t)\,
\biggl[\frac{1}{t}-\frac{\sin Nt}{Nt^2}\biggr]\,dt.\end{eqnarray*}
Treating $J_1$ as $I_1$ above, we get

\begin{eqnarray*}J_1&=&-\frac{1}{\pi}\lim\limits_{N\to\infty}\int_{-\infty}^0 F_s(-v)\,
\biggl[\frac{1}{v-x}-\frac{\sin(v-x)N}{N(v-x)^2}\biggr]\,dv\\
&=&\frac{1}{\pi}\lim\limits_{N\to\infty}\int_{-\infty}^0 F_s(v)\,
\biggl[\frac{1}{v-x}-\frac{\sin(v-x)N}{N(v-x)^2}\biggr]\,dv.\end{eqnarray*}
The last equality follows from the oddness of $F_s.$  Substituting
$v-x=t,$ we obtain

\begin{eqnarray*}J_1=\frac{1}{\pi}\lim\limits_{N\to\infty}\int_{-\infty}^{-x} F_s(x+t)\,
\biggl[\frac{1}{t}-\frac{\sin Nt}{Nt^2}\biggr]\,dt.\end{eqnarray*}
Therefore,

\begin{eqnarray*}F_c(x)=\frac{1}{\pi}\lim\limits_{N\to\infty}\int_{-\infty}^\infty F_s(x+t)
\biggl[\frac{1}{t}-\frac{\sin Nt}{Nt^2}\biggr]\,dt.\end{eqnarray*}
Finally, it follows from Theorem A that for almost all $x$)

\begin{eqnarray}\label{sico}F_c(x)=-{\mathcal H}F_s(x).\end{eqnarray}
This completes the proof. \hfill$\Box$

\medskip

Let us comment on the obtained results. In fact, the proof of
Theorem \ref{ans} shows that more general results than stated are
obtained. Indeed, formulas (\ref{cosi}) and (\ref{sico}) are more
informative than the assertion of Theorem \ref{ans}. To be precise,
such formulas are known, see \cite[(5.42) and (5.43)]{King}.
However, the situation is much more delicate. These formulas are
proved in \cite{King} for square integrable functions by applying
the Riemann-Lebesgue lemma in an appropriate place (5.44). But in
\cite[\S 6.19]{King} more details are given (see also \cite{Dyn})
and it is shown that the possibility to apply the Riemann-Lebesgue
lemma in that argument is equivalent to (Carleson's solution of)
Lusin's conjecture. In our $L^1$ setting this is by no means
applicable. And, indeed, our proof is different and rests on less
restrictive Theorem A. This is well agrees with what E.M. Dyn'kin
wrote in his well-known survey on singular integrals \cite{Dyn}: "In
fact, the theory of singular integrals is a technical subject where
ideas cannot be separated from the techniques."

\bigskip

\section{Sufficient conditions}

Analyzing the proof in \cite{IT}, one can see that in fact their
results are similar to ours, that is, can also be given in terms of
the (discrete) Hilbert transform. In that case (\ref{condser}) is
simply a sufficient condition for the summability of the discrete
Hilbert transform.

An analog of (\ref{condser}) for functions cannot be the only
sufficient condition for the integrability of the Hilbert transform.
Indeed, a known counter-example of the indicator function of an
interval works here as well: of course, it stands up to the
multiplication by logarithm.

Thus, we are going to give sufficient conditions for the
integrability of the Hilbert transforms in the spirit of those in
\cite{IT} for sequences supplied by some additional properties.

\subsection{General conditions}

First of all, examining integrability of the Hilbert transform, one
can test the integral over, say, $|t|\le\frac{3}{2}|x|.$ Indeed, for
$x>0,$ we have

\begin{eqnarray*}&\quad&\int_0^\infty\,\biggl|\,\int_{\frac32x}^\infty\frac{g(t)}{x-t}\,dt\,\biggr|\,dx\\
&\le&\int_0^\infty\,|g(t)|\,\int_0^{2t/3}\,\frac{dx}{x-t}\,dt=\ln3\int_0^\infty\,|g(t)|\,dt.
\end{eqnarray*}
The rest is estimated in a similar manner.

Further, since

\begin{eqnarray}\label{canc}\int_{-a}^{a}\,\frac{dt}{x-t}=0\end{eqnarray}
for any $a>0$ when the integral is understood in the principal value
sense, we can always consider

\begin{eqnarray*}\int_{-a}^{a}\,\frac{g(t)-g(x)}{x-t}\,dt\end{eqnarray*}
instead of the Hilbert transform truncated to $[-a,a].$

When in the definition of the Hilbert transform (\ref{dht}) the
function $g$ is odd, we will denote this transform by
$\mathcal{H}_o,$ and it is equal to

\begin{eqnarray}\label{oht}{\mathcal{H}_o}g(x)=\frac{2}{\pi}
\int_0^\infty\frac{tg(t)}{x^2-t^2}\,dt.           \end{eqnarray}
When $g$ is even its Hilbert transform ${\mathcal H}_e$ can be
rewritten as

\begin{eqnarray}\label{eht}{\mathcal{H}_e}g(x)=\frac{2}{\pi}
\int_0^\infty\frac{xg(t)}{x^2-t^2}\,dt.             \end{eqnarray}
Of course, both integrals should be understood in the principal
value sense (see, e.g., \cite[Ch.4, \S 4.2]{King}).

Since the functions $F_c$ and $F_s$ are even and odd, respectively,
their Hilbert transforms can be represented as

\begin{eqnarray}\label{Fteht}{\mathcal{H}_e}F_c(x)=\frac{2}{\pi}
\int_0^\infty\frac{xF_c(t)}{x^2-t^2}\,dt.\end{eqnarray} and

\begin{eqnarray}\label{Ftoht}\mathcal{H}_oF_s(x)=\frac{2}{\pi}
\int_0^\infty\frac{tF_s(t)}{x^2-t^2}\,dt.         \end{eqnarray}
These may be useful in some applications.

In fact, the most known condition is the following. If $g$ is of
compact support, a classical Zygmund $L\log L$ condition (see, e.g.,
\cite{Zg}) ensures the integrability of the Hilbert transform. More
precisely, the condition is the integrability of $g\log^+|g|,$ where
the $\log^+|g|$ notation means $\log|g|$ when $|g|>1$ and $0$
otherwise. As E.M. Stein has shown in \cite{St69}, this condition is
necessary on the intervals where the function is positive.

However, this condition looks quite restrictive in our case. We will
prove the following result.

\begin{theorem}\label{sciht} Let $g$ be an integrable function on
$\mathbb R$ which satisfies conditions (\ref{vm}),

\begin{eqnarray}\label{logc}\int_{|x|\ge1/2}|g(x)|\log3|x|\,dx\end{eqnarray}
and

\begin{eqnarray}\label{local}\int_{\mathbb
R}\,\int_{-\frac12\min(|x|,1)}^{\frac12\min(|x|,1)}\,\biggl|
\frac{g(x+t)-g(x)}{t}\biggr|\,dt\,dx.             \end{eqnarray}
Then $g\in H^1(\mathbb R)$.
\end{theorem}

Since each function can be represented as the sum of its even and
odd parts, we will prove Theorem \ref{sciht} separately for odd and
even functions. Thus, from now on we can consider $g$ to be defined
on $\mathbb R_+$ and analyze either (\ref{oht}) or (\ref{eht})
rather than the general Hilbert transform.

\subsection{Odd functions}

Though an odd function always satisfies (\ref{vm}), not every odd
integrable function belongs to $H^1(\mathbb R)$, for a
counterexample see, e.g., \cite{LiTi0}. Paley-Wiener's theorem (see
\cite{PW}; for alternative proof and discussion, see Zygmund's paper
\cite{Zyg}) asserts that if $g\in L^1(\mathbb R)$ is an odd and
monotone decreasing on $\mathbb{R}_+$ function, then
$\mathcal{H}g\in L^1$, i.e., $g$ is in $H^1(\mathbb R)$. Recently,
in \cite[Thm.6.1]{LiTi1}, this theorem has been extended to a class
of functions more general than monotone ones. However, it is
doubtful that these results are really practical in our situation.

Back to Theorem \ref{sciht}, we can consider

\begin{eqnarray*}\frac{2}{\pi}\int_{x/2}^{3x/2}\frac{tg(t)}{x^2-t^2}\,dt
\end{eqnarray*}
instead of (\ref{oht}). Indeed, the possibility of restricting to
that upper limit has been justified above. Similarly,

\begin{eqnarray}\label{avsm}\int_0^\infty\,\biggl|\int_0^{x/2}\frac{tg(t)}{x^2-t^2}\,dt
\biggr|\,dx&\le&\int_0^\infty|g(t)|t\int_{2t}^\infty\frac{dx}{x^2-t^2}\,dt\nonumber\\
&\le&\frac23\int_0^\infty|g(t)|\,dt.                  \end{eqnarray}
Now, like in (\ref{canc}), we have

\begin{eqnarray}\label{small}\int_0^1\,\biggl|\int_{x/2}^{3x/2}\frac{tg(t)}{x^2-t^2}\,dt
\biggr|\,dx&\le&\int_0^1\,\biggl|\int_{x/2}^{3x/2}\frac{t[g(t)-g(x)]}{x^2-t^2}\,dt
\biggr|\,dx+O(\int_0^\infty|g(t)|\,dt)\nonumber\\
&\le&\int_0^1\,\int_{-x/2}^{x/2}\frac{|g(x+t)-g(x)|}{|t|}\,dt\,dx
+O(\int_0^\infty|g(t)|\,dt).                         \end{eqnarray}
When $x\ge1$ we first estimate

\begin{eqnarray*}\int_1^\infty\,\biggl|\int^{x-1/2}_{x/2}\frac{tg(t)}{x^2-t^2}\,dt
\biggr|\,dx&\le&\int_{1/2}^\infty\,t|g(t)|\,\int^{2t}_{t+1/2}\frac{dx}{x^2-t^2}\,dt\\
&\le& C \int_{1/2}^\infty |g(t)|\ln3t\,dt\end{eqnarray*}      and

\begin{eqnarray*}\int_1^\infty\,\biggl|\int_{x+1/2}^{3x/2}\frac{tg(t)}{x^2-t^2}\,dt
\biggr|\,dx&\le&\int_{3/2}^\infty\,t|g(t)|\,\int_{2t/3}^{t-1/2}\frac{dx}{t^2-x^2}\,dt\\
&\le& C \int_{1/2}^\infty |g(t)|\ln3t\,dt.       \end{eqnarray*}
These two bounds lead to the logarithmic condition (\ref{logc}). The
remained integral

\begin{eqnarray*}\int_1^\infty\,\biggl|\int^{x-1/2}_{x+1/2}\frac{tg(t)}{x^2-t^2}\,dt\biggr|\,dx\end{eqnarray*}
is estimated exactly like that in (\ref{small}). Applying
(\ref{canc}), we obtain

\begin{eqnarray}\label{big}\int_1^\infty\,\biggl|\int_{x-1/2}^{x+1/2}\frac{tg(t)}{x^2-t^2}\,dt
\biggr|\,dx&\le&\int_1^\infty\,\biggl|\int_{x-1/2}^{x+1/2}\frac{t[g(t)-g(x)]}{x^2-t^2}\,dt
\biggr|\,dx+O(\int_0^\infty|g(t)|\,dt)\nonumber\\
&\le&\int_1^\infty\,\int_{-1/2}^{1/2}\frac{|g(x+t)-g(x)|}{|t|}\,dt\,dx+O(\int_0^\infty|g(t)|\,dt).\end{eqnarray}
Combining all the obtained estimates, we arrive at the required
result.

\subsection{Even functions}

While an odd function always satisfies (\ref{vm}), in the case of
even functions the situation is more delicate: the function must
satisfy (\ref{vm}) already on the half-axis. With this in hand, the
proof goes along the same lines as that for odd functions. The only
problem is that an estimate like (\ref{avsm}) does not follow
immediately from the formula (\ref{eht}). However, using the above
remark on the cancelation property for $g$ on the half-axis, we can
rewrite (\ref{eht}) as

\begin{eqnarray}\label{eht1}{\mathcal{H}_e}g(x)=\frac{2}{\pi}
\int_0^\infty\,g(t)\biggl[\frac{x}{x^2-t^2}-\frac{2}{x}\biggr]\,dt.
\end{eqnarray}
Now,

\begin{eqnarray}\label{avsm1}\int_0^\infty\,\biggl|\int_0^{x/2}\frac{t^2g(t)}{x(x^2-t^2)}\,dt
\biggr|\,dx&\le&\int_0^\infty|g(t)|t^2\int_{2t}^\infty\frac{dx}{x(x^2-t^2)}\,dt\nonumber\\
&\le&\frac16\int_0^\infty|g(t)|\,dt.\end{eqnarray}

This additional term $\frac2x$ does not affect the other estimates
of the previous subsection.

The proof is complete.

\bigskip

\section{Concluding remarks}

First of all, the relations (\ref{cosi}) and (\ref{sico}) are of
interest by their own.

The assertions of Theorem \ref{ans} can be reformulated in terms of
Hardy spaces: {\it belonging of $F_c$ ($F_s)$ to the real Hardy
space $H^1(\mathbb R)$ ensures the integrability of} $F_s$ ($F_c$).

The problem of sharpness is simple in this case: any known
counterexample of an integrable function with non-integrable Hilbert
transform works perfectly. For example, let
$F_c(x)=\frac{1}{1+x^2},$ the Fourier transform of $f(t)=e^{-|t|}.$
Surely, $F_c\in L^1(\mathbb R).$ However, $f$ cannot be re-expanded
in the integrable sine Fourier transform, since ${\mathcal H}F_c(x)=
\frac{x}{1+x^2}\not\in L^1(\mathbb R).$ That this is true, one can
see from the fact that the odd extension of this $F_c$ from the
right half-axis to the whole $\mathbb R$ is not continuous at zero.

More can be said about odd functions. Certain convenient conditions
for belonging of such functions to $H^1(\mathbb R)$ are known for
quite a long time. They are functions (Fourier transform) analogs of
important sufficient sequence conditions for the integrability of
trigonometric series (see, e.g., \cite{Te1} and \cite{Fomin}) and
can be found, for example, in \cite{L0} and in \cite{GM1}. In fact,
many of these subspaces first appeared in \cite{Borw}. For
$1<q\le\infty,$ set

\begin{eqnarray*}\|g\|_{A_q}=\int_0^\infty\left(\frac{1}{u}\int_{u\le|t|\le 2u}
|g(t)|^qdt \right)^{1/q}du,                       \end{eqnarray*}
with a standard modification when $q=\infty.$ In other words,
belonging of $g$ to one of the spaces $A_q$ ensures the
integrability of the odd Hilbert transform of $g.$

\bigskip

\section{Acknowledgements}

The author's attention to this problem was brought by R.M. Trigub.
Thanks to him and to J.-P. Kahane for stimulating discussions.

\end{document}